\documentclass[11pt]{article}
\usepackage{amsmath}
\usepackage{amssymb}
\usepackage{amsthm}
\usepackage{url}

\newcommand\sa{\smallskipamount}
\newcommand\ma{\medskipamount}
\newcommand\ba{\bigskipamount}
\newcommand\sLP{\\[\sa]}
\newcommand\mLP{\\[\ma]}
\newcommand\bLP{\\[\ba]}
\newcommand\bPP{\\[\ba]\indent}
\newcommand\CC{\mathbb{C}}
\newcommand\ZZ{\mathbb{Z}}
\newcommand\FSJ{{\cal J}}
\newcommand\FSL{{\cal L}}
\newcommand\al\alpha
\newcommand\be\beta
\newcommand\de\delta
\newcommand\tha\theta
\newcommand\la\lambda

\newcommand\half{\frac12}
\newcommand\thalf{\tfrac12}
\newcommand\iy\infty
\newcommand\pa\partial
\newcommand\Znonneg{\ZZ_{\ge0}}
\newcommand{\hyp}[5]{\,\mbox{}_{#1}F_{#2}\!\left(
  \genfrac{}{}{0pt}{}{#3}{#4};#5\right)}
\newcommand{\qhyp}[5]{\,\mbox{}_{#1}\phi_{#2}\!\left(
  \genfrac{}{}{0pt}{}{#3}{#4};#5\right)}
\newcommand\LHS{left-hand side}
\newcommand\RHS{right-hand side}

\numberwithin{equation}{section}
\newtheorem{theorem}{Theorem}[section]
\newtheorem{proposition}[theorem]{Proposition}
\newtheorem{lemma}[theorem]{Lemma}
\newtheorem{Definition}[theorem]{Definition}

\newtheorem{Remark}[theorem]{Remark}
\newenvironment{remark}{\begin{Remark}\rm}{\end{Remark}}
\newcommand\Proof{\noindent{\bf Proof}\quad}

\begin{document}
\title{Explicit matrix inverses for lower triangular matrices with entries
involving Jacobi polynomials}
\author{Leandro Cagliero\footnote{Partially supported by grants from CONICET and SECyT-UNC.}  $\;$and Tom H. Koornwinder}
\date{\em Dedicated to Richard Askey on the occasion of his 80th birthday}
\maketitle
\begin{abstract}
For a two-parameter family of
lower triangular matrices with entries involving
Jacobi polynomials an explicit inverse is given, with entries involving
a sum of two Jacobi polynomials.
The formula simplifies in the Gegenbauer case and then
one choice of the parameter solves an open problem
in a recent paper by Koelink, van Pruijssen \& Rom\'an.
The two-parameter family is closely related to two two-parameter groups
of lower triangular matrices, of which we also give the explicit generators.
Another family of pairs of mutually inverse lower triangular matrices
with entries involving Jacobi polynomials, unrelated to the family just
mentioned, was given by
J.~Koekoek \& R.~Koekoek (1999). We show that this last family
is a limit case of a pair of connection relations between Askey-Wilson
polynomials having one of their four parameters in common.
\end{abstract}
%
%
%
\section{Introduction}
This note started as a kind of supplement
to the paper \cite{2} by
Koelink, van Pruijssen \& Rom\'an, but gradually it got a wider scope.
As for \cite{2} it solves an open problem there
(see Theorem~2.1 and paragraph after Theorem~6.2 in \cite{2})
to invert a
lower triangular matrix with entries involving Gegenbauer polynomials.
For a two-parameter family of such matrices involving Jacobi polynomials
we give the explicit inverse matrix in Theorem \ref{18}.
Specialization to Gegenbauer polynomials then gives a one-parameter family.
One specialization of the parameter in the latter family
gives the inversion desired in \cite{2}. Another specialization gives a matrix
inversion already handled by Brega \& Cagliero \cite{1}.

Our two-parameter family of Jacobi polynomials is closely related to
two commutative two-parameter groups of lower triangular matrices involving
Jacobi polynomials. We also give the explicit infinitesimal
generators of these two-parameter
groups. Furthermore we obtain a biorthogonality relation
for two explicit systems of functions on $\ZZ$ involving Jacobi polynomials
with respect to an explicit bilinear form on $\ZZ$.

Another two-parameter family of pairs of mutually inverse lower
triangular matrices
with entries involving Gegenbauer polynomials, unrelated to the family
mentioned above, is implied by Brown \& Roman \cite[(4.14)]{3}.
J.~Koekoek \& R.~Koekoek \cite[(17)]{9}, unaware of \cite{3},
generalized a one-parameter subfamily of this two-parameter family
to entries involving
Jacobi polynomials. We will show that this last family can be realized
as a limit case of a pair of connection relations between Askey-Wilson
polynomials having one of their four parameters in common. These Askey-Wilson
connection coefficients were first given by Askey \& Wilson
\cite[(6.5)]{5}. The limit
case connects Jacobi polynomials $P_n^{(\al,\be)}$ with shifted monomials
$x\mapsto (x-y)^k$.

The contents of the paper are as follows.
In Section \ref{59} some preliminaries about Jacobi polynomials are given.
Degenerate cases of Jacobi polynomials are classified in Section \ref{80}.
The main results about the mutually inverse lower triangular matrices
are stated in Section \ref{19}. This section ends with some open
problems.
The computations leading to the explicit inverse matrix of the first
family of lower triangular matrices are given in Section \ref{16}.
The two-parameter groups and their generators are treated in Section \ref{60}.
The biorthogonal systems with respect to an explicit bilinear form are
the topic of Section \ref{81}.
Finally,
the computations giving the limit of the Askey-Wilson connection relations are
done in Section \ref{35}.

The reader may start in Section \ref{19} and then continue with
Section \ref{16} or with Sections \ref{60} and~\ref{81} or
with Section \ref{35}.
The preliminary sections \ref{59} and \ref{80} can be consulted when needed.
\mLP
{\em Acknowledgements}\quad
We thank the referees for careful reading and in particular one referee for suggesting
Remark \ref{89}.

We thank Michael Schlosser for the suggestion to look for a limit case of
the Askey-Wilson connection relations, and we thank Roelof Koekoek for
calling our attention to \cite{9}.
\section{Preliminaries about Jacobi polynomials}
\label{59}
Jacobi polynomials (see for instance \cite[Chapter IV]{12},
\cite[Chapter 6]{15}, \cite[Chapter 4]{7}, \cite[Section 9.8]{4},
\cite[Chapter 18]{14})
can be expressed in terms of the Gauss hypergeometric function
by
\begin{multline}
P_n^{(\al,\be)}(x):=\frac{(\al+1)_n}{n!}\,
\hyp21{-n,n+\al+\be+1}{\al+1}{\thalf(1-x)}\\
=\sum_{k=0}^n\frac{(n+\al+\be+1)_k\,\,(\al+k+1)_{n-k}}{k!\,(n-k)!}\,
\left(\frac{x-1}2\right)^k.
\label{10}
\end{multline}
Note that they are well-defined for all values of $\al,\be$. Their
normalization avoids
artificial singularities. Jacobi polynomials satisfy a Rodrigues formula
\begin{equation}
P_n^{(\al,\be)}(x)=\frac{(-1)^n}{2^n n!}\,(1-x)^{-\al}\,(1+x)^{-\be}\,
\left(\frac d{dx}\right)^n\,
\Bigl((1-x)^{n+\al}(1+x)^{n+\be}\Bigr).
\label{11}
\end{equation}

For $\al=\be$ Jacobi polynomials are often
written as Gegenbauer polynomials:
\begin{multline}
C_n^{(\la)}(x):=\frac{(2\la)_n}{(\la+\thalf)_n}\,P_n^{(\la-\half,\la-\half)}(x)
=\sum_{k=0}^{[n/2]}\frac{(-1)^k(\la)_{n-k}}{k!\,(n-2k)!}\,(2x)^{n-2k}\\
=\frac{2^n(\la)_n}{n!}\,x^n\,
\hyp21{-\thalf n,-\thalf n+\thalf}{1-\la-n}{\frac1{x^2}},
\label{48}
\end{multline}
where we also used \cite[10.9(18)]{13}.
Thus $C_n^{(0)}(x)=\de_{n,0}$, which will be kept as a convention in
this paper,
although in the literature the case $\la=0$ is usually rescaled in order to obtain
the Chebyshev polynomials of the first kind.
In the proportionality factor in the second part of \eqref{48}
artificial singularities
can occur. This factor should be understood by continuity in $\la$. We can
rewrite the first equality in \eqref{48} as
\begin{equation*}
C_{2m}^{(\la)}(x)=\frac{2^{2m}(\la)_m}{(\la+m+\thalf)_m}\,
P_{2m}^{(\la-\half,\la-\half)}(x),\quad
C_{2m-1}^{(\la)}(x)=\frac{2^{2m-1}(\la)_m}{(\la+m-\thalf)_m}\,
P_{2m-1}^{(\la-\half,\la-\half)}(x)
\end{equation*}
or as
\begin{equation}
C_n^{(\la)}(x)=\frac{2^{2n}(\la)_n}{(n+2\la)_n}\,
P_n^{(\la-\half,\la-\half)}(x).
\label{67}
\end{equation}

In the Legendre case $\al=\be=0$ we write
$P_n(x):=P_n^{(0,0)}(x)$.
There are symmetries
\begin{equation}
P_n^{(\al,\be)}(x)=(-1)^n P_n^{(\be,\al)}(-x),\qquad
C_n^{(\la)}(x)=(-1)^n C_n^{(\la)}(-x).
\label{49}
\end{equation}

For Jacobi polynomials we will need the following generating function
(see \cite[(4.4.5)]{12}):
\begin{equation}
\sum_{n=0}^\iy P_n^{(\al,\be)}(x)\,w^n
=2^{\al+\be}R^{-1}(1-w+R)^{-\al}(1+w+R)^{-\be},
\quad R:=(1-2xw+w^2)^\half,
\label{50}
\end{equation}
convergent for $x\in[-1,1]$, $|w|<1$.
A more simple generating function for Gegenbauer polynomials (but not the
case $\al=\be$ of \eqref{50}) is the following (see \cite[(4.7.23)]{12}):
\begin{equation}
\sum_{n=0}^\iy C_n^{(\la)}(x)\,w^n=(1-2xw+w^2)^{-\la}\qquad
(x\in[-1,1],\;|w|<1).
\label{51}
\end{equation}
\section{Degenerate cases of Jacobi polynomials}
\label{80}
This section is not needed very much in the sequel. It may be skipped on first
reading.

For $\al,\be>-1$ Jacobi polynomials are orthogonal on the interval $(-1,1)$
with respect to the weight function $(1-x)^\al(1+x)^\be$, but we will not
deal with this property in the paper. However, since in our formulas
$\al,\be$ will be allowed to be arbitrarily complex, and definitely not
only larger than $-1$, it is relevant to see which degeneracies can
occur in \eqref{10}, i.e., when coefficients in the sum
on the right of \eqref{10}
become zero (here we assume $n>0$).
There are two shifted factorials in the numerator of the terms which
can cause this:
\begin{enumerate}
\item
$(n+\al+\be+1)_k=0$ for some $k\in\{1,\ldots,n\}$, i.e.,
$(n+\al+\be+1)_n=0$, i.e,
$\al+\be\in\ZZ_{\le-2}$, $n+\al+\be+1\le0$ and $2n+\al+\be\ge0$.
Then $(n+\al+\be+1)_k=0$ for $k=-n-\al-\be,\ldots,n$.
\item
$(\al+k+1)_{n-k}=0$ for some $k\in\{0,1,\ldots,n-1\}$,
i.e., $(\al+1)_n=0$, i.e.,
$\al\in\ZZ_{\le-1}$ and $n+\al\ge0$. Then $(\al+k+1)_{n-k}=0$ for
$k=0,\ldots,-\al-1$.
\end{enumerate}
By combining these two cases we see when $(n+\al+\be+1)_k(\al+k+1)_{n-k}=0$
for all $k\in\{0,\ldots,n\}$:

\begin{proposition}
\label{65}
$P_n^{(\al,\be)}(x)=0$ identically in $x$ iff $\al,\be\in\ZZ_{\le-1}$ and
$\max(-\al,-\be)\le n\le-\al-\be-1$.
\end{proposition}

Case 1 above causes that $P_n^{(\al,\be)}(x)$ has degree lower than $n$
in $x$,
while case 2 causes that $P_n^{(\al,\be)}(x)$ vanishes for $x=1$
with a certain
multiplicity. A similar case with vanishing at $-1$ then follows by \eqref{49}.
In all these cases we can look at the \RHS\ of \eqref{10} in a different way
and thus obtain a transformation formula such that the true degree or the
multiplicity of vanishing at 1 or $-1$ can be read off from the transformed
expression. The results are:

\begin{proposition}
\label{66}
Let $n>0$. Assume that $P_n^{(\al,\be)}(x)$ does not vanish identically in $x$.
\sLP
{\bf(a)} $P_n^{(\al,\be)}(x)$ has degree $<n$ in $x$ iff
$\al+\be\in\ZZ_{\le-2}$,
$n+\al+\be+1\le0$ and $2n+\al+\be\ge0$. Then the degree is
$-n-\al-\be-1$ and
\begin{equation*}
P_n^{(\al,\be)}(x)=
\frac{(-n-\be)_{2n+\al+\be+1}}{(-n-\al-\be)_{2n+\al+\be+1}}
\,P_{-n-\al-\be-1}^{(\al,\be)}(x).
\end{equation*}
{\bf(b)} $P_n^{(\al,\be)}(1)=0$ iff $\al\in\ZZ_{\le-1}$ and $n+\al\ge0$.
Then the zero
at 1 has multiplicity $-\al$ and
\begin{equation*}
P_n^{(\al,\be)}(x)=
\frac{(n+\al+\be+1)_{-\al}\,(n+\al)!}{n!}\,
\left(\frac{x-1}2\right)^{-\al}
\,P_{n+\al}^{(-\al,\be)}(x).
\end{equation*}
{\bf(c)} $P_n^{(\al,\be)}(-1)=0$ iff $\be\in\ZZ_{\le-1}$ and $n+\be\ge0$.
Then the zero at $-1$ has multiplicity $-\be$ and
\begin{equation*}
P_n^{(\al,\be)}(x)=
\frac{(n+\al+\be+1)_{-\be}\,(n+\be)!}{n!}\,
\left(\frac{x+1}2\right)^{-\be}
P_{n+\be}^{(\al,-\be)}(x).
\end{equation*}
\end{proposition}

Combinations of the cases in this last proposition can occur. Then the
corresponding
transformation formulas can be combined. For instance, combination of
(a) and (b) yields:
\bLP
{\bf(d)} $\al,\be\in\ZZ$ and $\be+2\le\al\le-1$. Then for
$\max(-\al,-\thalf(\al+\be))\le n\le-\be-1$ we have
\begin{equation*}
P_n^{(\al,\be)}(x)=\left(\frac{1-x}2\right)^{-\al}
P_{-n-\be-1}^{(-\al,\be)}(x).
\end{equation*}
A further combination of (d) with (c) is empty.
The combination of (a) and (c) can be obtained from (d) by using \eqref{49}:
\bLP
{\bf (e)} $\al,\be\in\ZZ$ and $\al+2\le\be\le-1$. Then for
$\max(-\be,-\thalf(\al+\be))\le n\le-\al-1$ we have
\begin{equation*}
P_n^{(\al,\be)}(x)=(-1)^{\al+1}\left(\frac{1+x}2\right)^{-\be}
P_{-n-\al-1}^{(\al,-\be)}(x).
\end{equation*}
Combination of (b) and (c) yields:
\bLP
{\bf (f)} $\al,\be\in\ZZ_{\le-1}$. Then for $n\ge-\al-\be$ we have
\begin{equation}
P_n^{(\al,\be)}(x)=\left(\frac{x-1}2\right)^{-\al}
\left(\frac{x+1}2\right)^{-\be}P_{n+\al+\be}^{(-\al,-\be)}(x).
\label{68}
\end{equation}

We can also consider vanishing of coefficients in the sum in \eqref{48}.
Let us rewrite this as a summation formula for Jacobi polynomials
$P_n^{(\al,\al)}(x)$
and let us distinguish between cases $n=2m$ and $n=2m-1$:
\begin{align}
P_{2m}^{(\al,\al)}(x)&=2^{-2m}(\al+m+1)_m\,
\sum_{k=0}^m\frac{(-1)^k(\al+m+\thalf)_{m-k}}{k!\,(2m-2k)!}\,(2x)^{2m-2k},
\label{63}\\
P_{2m-1}^{(\al,\al)}(x)&=2^{1-2m}(\al+m)_m\,
\sum_{k=0}^{m-1}\frac{(-1)^k(\al+m+\thalf)_{m-1-k}}{k!\,(2m-1-2k)!}\,
(2x)^{2m-1-2k}.
\label{64}
\end{align}
From \eqref{63} and \eqref{64} Proposition \ref{65} and Proposition \ref{66}(a)
can again be derived in the
case $\al=\be$. Furthermore, we conclude that, if \eqref{63} and \eqref{64} are
not identically zero in $x$, then they have no zero at $x=0$ (in case of
\eqref{63}) respectively no zero of multiplicity higher than one at $x=0$
(in case
of \eqref{64}).
\section{Main results}
\label{19}
In Section \ref{16} it will be shown that
\begin{equation}
\sum_{k=0}^n P_{n-k}^{(\al_1+k,\be_1+k)}(x)\,
P_k^{(\al_2-k,\be_2-k)}(x)
=P_n^{(\al_1+\al_2,\be_1+\be_2)}(x)
\label{12}
\end{equation}
and
\begin{multline}
\sum_{k=0}^n k\,P_{n-k}^{(\al_1+k,\be_1+k)}(x)\,
P_k^{(\al_2-k,\be_2-k)}(x)
=\frac{n(\al_2+\be_2)}{\al_1+\al_2+\be_1+\be_2+2n}\,
P_n^{(\al_1+\al_2,\be_1+\be_2)}(x)\\
+\frac{\al_2\be_1-\al_1\be_2+n(\al_2-\be_2)}{\al_1+\al_2+\be_1+\be_2+2n}\,
P_{n-1}^{(\al_1+\al_2,\be_1+\be_2)}(x)\quad(n>0).
\label{14}
\end{multline}
For $\al_1=-\al_2=\al$, $\be_1=-\be_2=\be$ these formulas reduce to
\begin{align}
\sum_{k=0}^n P_{n-k}^{(\al+k,\be+k)}(x)\,
P_k^{(-\al-k,-\be-k)}(x)&=P_n(x),
\label{40}\\
\sum_{k=0}^n k\,P_{n-k}^{(\al+k,\be+k)}(x)\,P_k^{(-\al-k,-\be-k)}(x)
&=-\thalf(\al+\be)P_n(x)+\thalf(\be-\al)P_{n-1}(x)\quad(n>0).
\label{41}
\end{align}
From \eqref{40}, \eqref{41} and \eqref{12} we obtain for $n>0$ that
\begin{multline}
\sum_{k=0}^n P_{n-k}^{(\al+k,\be+k)}(x)\big((k+\be)
P_k^{(-\al-k,-\be-k)}(x)+(\al-\be)P_k^{(-\al-k,-\be-k-1)}(x)\big)\\
=-\thalf(\al+\be)P_n(x)+\thalf(\be-\al)P_{n-1}(x)
+\be P_n(x)+(\al-\be)P_n^{(0,-1)}(x)=0
\label{86}
\end{multline}
by \cite[10.8(36)]{13}. Thus we have derived for $n=0,1,2,\ldots\;$ that
\begin{equation}
\sum_{k=0}^n P_{n-k}^{(\al+k,\be+k)}(x)\left(\frac{k+\be}\al\,
P_k^{(-\al-k,-\be-k)}(x)
+\frac{\al-\be}\al\,P_k^{(-\al-k,-\be-k-1)}(x)\right)=\de_{n,0}
\label{42}
\end{equation}
and, in particular,
\begin{equation}
\sum_{k=0}^n\frac{k+\al}\al\,
P_{n-k}^{(\al+k,\al+k)}(x)\,P_k^{(-\al-k,-\al-k)}(x)=
\de_{n,0}.
\label{5}
\end{equation}

Now make in \eqref{42} the substitutions $n\to m-n$, $k\to k-n$, $\al\to\al+n$,
$\be\to\be+n$, where the new variables $m,n$
can be arbitrarily integer such that $m\ge n$.
The resulting identity is:
\begin{equation}
\sum_{k=n}^m P_{m-k}^{(\al+k,\be+k)}(x)\left(
\frac{k+\be}{n+\al}\,
P_{k-n}^{(-\al-k,-\be-k)}(x)
+\frac{\al-\be}{n+\al}\,P_{k-n}^{(-\al-k,-\be-k-1)}(x)\right)=\de_{m,n}\quad
(m\ge n).
\label{33}
\end{equation}
In this and related formulas it turns out that the expression remains
continuous in
$\al$ as $\al$ tends to the apparent singularity, see Remark \ref{71}.

Let $\FSL_\iy$ be the group of all 
lower triangular $\iy\times\iy$ matrices
(doubly infinite, i.e., with indices running over all integers)
for which the entries depend on a complex
variable $x$ (usually polynomially), but which have the entries on the main
diagonal identically 1.
The identity \eqref{33} can be rephrased by giving two explicit elements of
$\FSL_\iy$ which are inverse to each other:
\begin{theorem}
\label{18}
 $LM=I=ML$ where $L=L^{(\al,\be)}$ and $M=M^{(\al,\be)}$
 are lower triangular matrices for which
 the lower triangular entries ($m\ge n$) are given by
\begin{equation}
L_{m,n}^{(\al,\be)}=P_{m-n}^{(\al+n,\be+n)}(x),\quad
M_{m,n}^{(\al,\be)}=\frac{m+\be}{n+\al}\,
P_{m-n}^{(-\al-m,-\be-m)}(x)
+\frac{\al-\be}{n+\al}\,P_{m-n}^{(-\al-m,-\be-m-1)}(x).
\label{20}
\end{equation}
\end{theorem}

In the Gegenbauer case $\al=\be$ formulas \eqref{33} and \eqref{20} simplify
because the term with factor $\al-\be$ vanishes.

The lower triangular matrices in Theorem \ref{18} can also be considered
with entries $m,n$ running over all integers $\ge n_0$ for some integer $n_0$,
in particular with entries running over all nonnegative integers.
There is no loss of generality in doing this because
\[
L_{m,n}^{(\al,\be)}=L_{m-k,n-k}^{(\al+k,\be+k)},\qquad
M_{m,n}^{(\al,\be)}=M_{m-k,n-k}^{(\al+k,\be+k)}.
\]

There are two different places in the literature where Theorem \ref{18}
can be used,
for $\al=\be=1$ and $\al=\be=\thalf$, respectively:
\begin{enumerate}
\item
The case $\al=\be=1$ (matrix entries running over all nonnegative integers)
occurs in Brega \& Cagliero \cite[p.471]{1}
with a proof similar as given here.
\item
The case $\al=\be=\thalf$ of the matrix $L^{(\al,\be)}$ in \eqref{20} occurs in
Koelink, van Pruijssen \& Rom\'an
\cite[Theorem 2.1]{2} in the form of the lower triangular matrix $L(x)$
given by
\begin{equation*}
(L(x))_{m,n}=\frac{m!\,(2n+1)!}{(m+n+1)!\,n!}\,C_{m-n}^{(n+1)}(x)
=\frac{m!\,(\tfrac32)_n}{(\tfrac32)_m\,n!}\,P_{m-n}^{(n+\half,n+\half)}(x)\quad
(m\ge n\ge0).
\end{equation*}
There the matrix has finite size (which does not matter for the
purpose of inversion).
As the authors wrote in \cite[paragraph after Theorem 6.2]{2},
they tried to find an explicit inverse matrix but
did not succeed. We can give the inverse by \eqref{20} for
$\al=\thalf$ as follows:
\begin{equation*}
(L^{-1}(x))_{m,n}=\frac{m!\,(m+n)!}{(2m)!\,n!}\,C_{m-n}^{(-m)}(x)=
\frac{m!\,(-\thalf)_{n+1}}{(-\thalf)_{m+1}\,n!}\,
P_{m-n}^{(-m-\half,-m-\half)}(x)
\quad(m\ge n\ge0).
\end{equation*}
Our result is mentioned in an Addendum at the end of \cite{2}.
\end{enumerate}
\begin{remark}
\label{71}
In the formula for $M_{m,n}$ in \eqref{20} the denominator only
gives an apparent singularity.
We have $M_{n,n}^{(\al,\be)}=1$ and for $m>n$ we obtain by
\eqref{10} that
\begin{multline*}
M_{m,n}^{(\al,\be)}=
-\sum_{k=0}^{m-n-1}\big((m+\be)(-\al-\be-m-n+1)_k
+(\al-\be)(-\al-\be-m-n)_k\big)\\
\times
\frac{(-\al-m+k+1)_{m-n-k-1}}{k!\,(m-n-k)!}\,\left(\frac{x-1}2\right)^k\\
-\,\frac{(\al+\be+2m)(-\al-\be-m-n+1)_{m-n-1}}{(m-n)!}\,
\left(\frac{x-1}2\right)^{m-n}.
\end{multline*}
\end{remark}
\begin{remark}
\label{89}
A referee suggested to see if \eqref{86}, which leads to Theorem \ref{18},
 can be generalized by starting in the left part with
\begin{equation}
\sum_{k=0}^n P_{n-k}^{(\al_1+k,\be_1+k)}(x)
\big((x_nk+y_n)P_k^{(\al_2-k,\be_2-k)}(x)+z_nP_k^{(\al_2-k,\be_2-k-1)}(x)\big)
\label{87}
\end{equation}
for yet unknown $x_n$, $y_n$, $z_n$ and then find for which choices of the unknowns
the generalization of the middle part of \eqref{86} will match with
\cite[10.8(36)]{13} and thus yield zero.
The solution (up to multiplication of the three unknowns by the same possibly
$n$-dependent factor) turns out to be 
\begin{equation}
x_n=-n-\alpha_1-\alpha_2, \quad
y_n=\beta_2\,n+\alpha_1\beta_2-\alpha_2\beta_1,\quad
z_n=(\alpha_2-\beta_2)n-\alpha_1\beta_2+\alpha_2\beta_1.
\label{88}
\end{equation}
Solutions independent of $n$ can be obtained iff $\al_1+\al_2=0=\be_1+\be_2$, precisely
the case which we already had in \eqref{86}. It is not clear how we can obtain a pair
of mutually inverse lower triangular
matrices in the general $n$-dependent case of \eqref{88} making
\eqref{87} zero.
\end{remark}
\begin{remark}
Just as we can go from  \eqref{42} to $LM=I$ and backwards, we can go
back and forth from $ML=I$ to the identity
\begin{equation*}
\sum_{k=0}^n P_k^{(\al,\be)}(x)\left(\frac{n+\be}{k+\al}\,
P_{n-k}^{(-\al-n,-\be-n)}(x)
+\frac{\al-\be}{k+\al}\,P_{n-k}^{(-\al-n,-\be-n-1)}(x)\right)=\de_{n,0}.
\end{equation*}
\end{remark}
\begin{remark}
\label{78}
By Theorem \ref{18} we have a biorthogonal system of functions on $\ZZ$ given
for $n\ge k$ by
\begin{align}
\phi_n^{(\al,\be,x)}(k)&:=L_{n,k}^{(\al,\be)}\;\;\,=P_{n-k}^{(\al+k,\be+k)}(x),
\label{75}\\
\psi_n^{(\al,\be,x)}(k)&:=M_{-k,-n}^{(\al,\be)}=
\frac{\be-k}{\al-n}\,
P_{n-k}^{(-\al+k,-\be+k)}(x)
+\frac{\al-\be}{\al-n}\,P_{n-k}^{(-\al+k,-\be+k-1)}(x),
\nonumber
\end{align}
and otherwise zero.
Here $\psi_n^{(\al,\be,x)}(k)$ for $\al\to n$ is to be understood as in
Remark \ref{71}.
Then
\begin{equation}
\sum_{k\in\ZZ}\phi_m^{(\al,\be,x)}(k)\,\psi_{-n}^{(\al,\be,x)}(-k)=\de_{m,n},
\label{69}
\end{equation}
where the sum actually only runs over $\{k\in\ZZ\mid n\le k\le m\}$.
If $\al$ and $\be$ are shifted by the same integer $j$ then the biorthogonal
system does not essentially change, since
\[
\phi_n^{(\al+j,\be+j,x)}(k)=\phi_{n+j}^{(\al,\be,x)}(k+j),\qquad
\psi_n^{(\al+j,\be+j,x)}(k)=\psi_{n-j}^{(\al,\be,x)}(k-j).
\]
For $\al=\be=0$ the biorthogonality relation \eqref{69} simplifies to
\begin{equation*}
\sum_{k\in\ZZ} \phi_m^{(0,0,x)}(k)\,
\phi_{-n}^{(0,0,x)}(-k)\,\frac kn=\de_{m,n}\,.
\end{equation*}
Again the sum only runs over $\{k\in\ZZ\mid n\le k\le m\}$
and the singularity for $n=0$ in $\phi_{-n}^{(0,0,x)}(-k)\,\frac kn$
is only apparent
because of Remark \ref{71}.
For $n\le 0\le -n<m$ the summation range is further restricted to
$\{k\in\ZZ\mid n\le k\le -n\}$ because of Proposition \ref{65}.
Similarly, for $n<-m\le 0\le m$ the summation range is restricted to
$\{k\in\ZZ\mid -m\le k\le m\}$.

From $ML=I$ we get a biorthogonality relation for the dual systems:
\begin{equation*}
\sum_{k\in\ZZ}\psi_{-k}^{(\al,\be,x)}(-m)\,\phi_k^{(\al,\be,x)}(n)
=\de_{m,n}\,.
\end{equation*}
\end{remark}
\begin{remark}
Brown \& Roman \cite[(4.12)--(4.15)]{3} obtain inverse relations involving
Gegenbauer polynomials of which a special case is close to \eqref{5}
but not equal to
it. It reads
\begin{equation}
\sum_{k=0}^n  \frac{(n+2\al+1)_k}{(2\al+2)_k}
P_{n-k}^{(\al+k,\al+k)}(x)\,P_k^{(-\al-k-1,-\al-k-1)}(x)=\de_{n,0}\,.
\label{30}
\end{equation}
In fact, they give a more general identity
\begin{equation}
\sum_{k=0}^n\frac\nu{\mu k+\nu}\,C_k^{(\mu k+\nu)}(x)\,
C_{n-k}^{(-\mu k-\nu)}(x)=\de_{n,0}\,.
\label{34}
\end{equation}
Then \eqref{30} is the case $\mu=-1$, $\nu=-\al-\thalf$ of \eqref{34},
while the case
$\mu=0$ of \eqref{34} is the very elegant formula
\begin{equation}
\sum_{k=0}^n C_k^{(\nu)}(x)\,
C_{n-k}^{(-\nu)}(x)=\de_{n,0}\,,
\label{43}
\end{equation}
which is also the case $\la=-\nu$  of the formula
\begin{equation}
\sum_{k=0}^n C_k^{(\nu)}(x)\,
C_{n-k}^{(\la)}(x)=C_n^{(\nu+\la)}(x),
\label{44}
\end{equation}
mentioned in \cite[(18.18.20)]{14}.
The identity \eqref{44} is a direct consequence
of the generating function~\eqref{51}.

Formula \eqref{30} is also the case $\al=\be$ of the identity
\begin{equation}
\sum_{k=0}^n \frac{(\al+\be+n+1)_k}{(\al+\be+2)_k}\,
P_{n-k}^{(\al+k,\be+k)}(x)\,
P_k^{(-\al-k-1,-\be-k-1)}(x)=\de_{n,0}.
\label{31}
\end{equation}
This last identity is a consequence of the pair of mutually inverse
lower triangular matrices \eqref{32}
implied by J. Koekoek \& R. Koekoek \cite[(17)]{9}.

In Section \ref{35} we will show that \eqref{31},
and hence \eqref{30}, is related
to a limit case of a connection formula for Askey-Wilson polynomials.
\end{remark}
\begin{remark}
There remain several interesting questions.
First of all, is there a larger family of explicit mutually inverse
lower triangular
matrices which includes both the family of Theorem \ref{18} and the family
\eqref{36} implying \eqref{31}? (See the attempt made in Remark \ref{89}.)$\;$
Furthermore, are there two simple systems
of special functions connected by the matrices in Theorem \ref{18}?
If yes, can this also be seen as a limit case for $q\to 1$ of some connection
formula in the $q$-case?
Concerning the pair of mutually inverse
lower triangular matrices \eqref{32} involving Jacobi polynomials
there are analogues for some other families of orthogonal polynomials
in the Askey scheme, for instance for Charlier and Meixner polynomials, as
surveyed by Koekoek \cite{10}. It would be interesting to see if these also
come from limit cases of the Askey-Wilson connection relations.
Finally there is the puzzling Brown-Roman formula
\eqref{34}.
Does this have an extension to Jacobi polynomials for general $\mu$?
It would also be interesting to generalize \eqref{44} such that it is
related to \eqref{34}.
\end{remark}
\section{Computations leading to Theorem \ref{18}}
\label{16}
\begin{lemma}
If the functions $f$ and $g$ have derivatives up to order $n$ then
\begin{align}
\sum_{k=0}^n\binom nk f^{(n-k)}(x)\,g^{(k)}(x)&=(fg)^{(n)}(x),
\label{8}\\
\sum_{k=0}^n k\binom nk f^{(n-k)}(x)\,g^{(k)}(x)&=n\,(fg')^{(n-1)}(x).
\label{9}
\end{align}
\end{lemma}
\Proof
Formula \eqref{8} is well-known. For the proof of \eqref{9}
rewrite its \LHS\ as
\[
n\,\sum_{j=0}^{n-1}\binom{n-1}j f^{(n-j-1)}(x)\,(g')^{(j)}(x)
\]
and use \eqref{8}.\qed
\bPP
By the Rodrigues formula \eqref{11} we have
\begin{multline*}
\sum_{k=0}^n P_{n-k}^{(\al_1+k,\be_1+k)}(x)\,
P_k^{(\al_2-k,\be_2-k)}(x)
=\frac{(-1)^n}{2^n n!}\,(1-x)^{-\al_1-\al_2}\,(1+x)^{-\be_1-\be_2}\\
\times \sum_{k=0}^n\binom nk
\left(\frac d{dx}\right)^{n-k}\,
\Bigl((1-x)^{n+\al_1}(1+x)^{n+\be_1}\Bigr)
\left(\frac d{dx}\right)^k\,
\Bigl((1-x)^{\al_2}(1+x)^{\be_2}\Bigr).
\end{multline*}
By \eqref{8} and again \eqref{11} we obtain \eqref{12}.

Similarly, by \eqref{11} and \eqref{9} we can write for $n>0$:
\begin{multline*}
\sum_{k=0}^n k\,P_{n-k}^{(\al_1+k,\be_1+k)}(x)\,
P_k^{(\al_2-k,\be_2-k)}(x)
=\frac{(-1)^n}{2^n n!}\,(1-x)^{-\al_1-\al_2}\,(1+x)^{-\be_1-\be_2}\\
\times \sum_{k=0}^n k\binom nk
\left(\frac d{dx}\right)^{n-k}\,
\Bigl((1-x)^{n+\al_1}(1+x)^{n+\be_1}\Bigr)
\left(\frac d{dx}\right)^k\,
\Bigl((1-x)^{\al_2}(1+x)^{\be_2}\Bigr)\\
=\frac{(-1)^n}{2^n (n-1)!}\,(1-x)^{-\al_1-\al_2}\,(1+x)^{-\be_1-\be_2}
\left(\frac d{dx}\right)^{n-1}\,
\Big((1-x)^{n+\al_1}(1+x)^{n+\be_1}\,
\tfrac d{dx}\big((1-x)^{\al_2}(1+x)^{\be_2}\big)\Big).
\end{multline*}
By straightforward computation we get
\begin{multline*}
(1-x)^{n+\al_1}(1+x)^{n+\be_1}\,
\frac d{dx}\Big((1-x)^{\al_2}(1+x)^{\be_2}\Big)\\
=\frac{\al_2+\be_2}{\al_1+\al_2+\be_1+\be_2+2n}\,
\frac d{dx}\Big((1-x)^{\al_1+\al_2+n}(1+x)^{\be_1+\be_2+n}\Big)\\
-2\,\frac{\al_2\be_1-\al_1\be_2+n(\al_2-\be_2)}{\al_1+\al_2+\be_1+\be_2+2n}\,
(1-x)^{\al_1+\al_2+n-1}(1+x)^{\be_1+\be_2+n-1}.
\end{multline*}
By \eqref{11} we finally obtain
\eqref{14}.
\section{Further matrix identities involving Jacobi polynomials}
\label{60}
As a consequence of the generating function \eqref{50} we have
\begin{equation*}
\sum_{n=0}^\iy\left(\sum_{k=0}^n
P_{n-k}^{(\al_1,\be_1)}(x)\,P_k^{(\al_2,\be_2)}(x)\right)w^n
=2^{\al_1+\al_2+\be_1+\be_2}R^{-2}(1-w+R)^{-\al_1-\al_2}(1+w+R)^{-\be_1-\be_2},
\end{equation*}
by which the inner sum on the \LHS\ as a function of $\al_1,\al_2,\be_1,\be_2$
only depends on $\al_1+\al_2,\be_1+\be_2$. In particular,
\begin{equation}
\sum_{k=0}^n P_{n-k}^{(\al_1,\be_1)}(x)\,P_k^{(\al_2,\be_2)}(x)=
\sum_{k=0}^n P_{n-k}(x)\,P_k^{(\al_1+\al_2,\be_1+\be_2)}(x).
\label{45}
\end{equation}
Formula \eqref{45} is quite similar to \eqref{12}. We can rewrite both
identities
as identities in $\FSL_\iy$ (the group of doubly infinite lower
triangular matrices
depending on
a complex variable $x$ and with 1 on the main diagonal).
Let $P^{(\al,\be)},Q^{(\al,\be)}\in\FSL_\iy$ with
\begin{equation}
P^{(\al,\be)}_{m,n}:=P_{m-n}^{(\al,\be)}(x),\quad
Q^{(\al,\be)}_{m,n}:=P_{m-n}^{(\al+n-m,\be+n-m)}(x)\qquad(m\ge n).
\label{52}
\end{equation}
Both are matrices of the form $A_{m,n}=f(m-n)$ (constant on each
diagonal, i.e.,
a Toeplitz matrix).
All such matrices in $\FSL_\iy$ commute.
Formulas \eqref{45} and \eqref{12} can be rephrased as:
\begin{equation}
P^{(\al_1,\be_1)}P^{(\al_2,\be_2)}=P^{(0,0)}P^{(\al_1+\al_2,\be_1+\be_2)},\quad
Q^{(\al_1,\be_1)}Q^{(\al_2,\be_2)}=Q^{(\al_1+\al_2,\be_1+\be_2)}.
\label{46}
\end{equation}
Also, by \eqref{10}, $Q^{(0,0)}=I$. Put
\begin{equation}
P_H^{(\al,\be)}:=P^{(\al,\be)}(P^{(0,0)})^{-1}=(P^{(0,0)})^{-1}P^{(\al,\be)},
\qquad Q_H^{(\al,\be)}:=Q^{(\al,\be)}.
\label{53}
\end{equation}
Then, by \eqref{46},
\begin{equation}
P_H^{(\al_1,\be_1)}P_H^{(\al_2,\be_2)}=P_H^{(\al_1+\al_2,\be_1+\be_2)},\qquad
Q_H^{(\al_1,\be_1)}Q_H^{(\al_2,\be_2)}=Q_H^{(\al_1+\al_2,\be_1+\be_2)}.
\label{47}
\end{equation}

Since $P_n(x)=C_n^{(\half)}(x)$ we see by \eqref{43} that
\begin{equation}
\big(\big(P^{(0,0)}\big)^{-1}\big)_{m,n}=C_{m-n}^{(-\half)}(x)\qquad(m\ge n).
\label{57}
\end{equation}
Also note that $C_0^{(-\half)}(x)=1$, $C_1^{(-\half)}(x)=-x$ and, for $n\ge2$,
\begin{align*}
C_n^{(-\half)}(x)&=\frac{2^n(-\thalf)_n}{n!}\,x^n\,
\hyp21{-\thalf n,\thalf-\thalf n}{\tfrac32-n}{\frac1{x^2}}\\
&=\frac{2^n(-\thalf)_n}{n!}\,(x^2-1)x^{n-2}\,
\hyp21{1-\thalf n,\tfrac32-\thalf n}{\tfrac32-n}{\frac1{x^2}}\\
&=\frac{1-x^2}{n(n-1)}\,C_{n-2}^{(\frac32)}(x)
=\frac{1-x^2}{2(n-1)}\,P_{n-2}^{(1,1)}(x).
\end{align*}
Here we used \eqref{48} and Euler's transformation formula
\cite[(2.2.7)]{15} for ${}_2F_1$ series.
Alternatively, use \eqref{67} and \eqref{68}.

By \eqref{47} the maps sending $(\al,\be)\in\CC^2$ to $P_H^{(\al,\be)}$ and to
$Q_H^{(\al,\be)}$ are both group homomorphisms from $\CC^2$ into
 $\FSL_\iy$.
The maps are entrywise analytic and entries on the \RHS s of \eqref{47}
are obtained from finite sums on the \LHS s. Thus we must have
\begin{equation}
P_H^{(\al,\be)}=\exp(\al A_P+\be B_P),\qquad
Q_H^{(\al,\be)}=\exp(\al A_Q+\be B_Q)
\label{55}
\end{equation}
for some strictly lower triangular matrices $A_P,B_P,A_Q,B_Q$,
and these matrices
can be computed by evaluating the derivatives
$\frac\pa{\pa\al}P_H^{(\al,\be)}$,
$\frac\pa{\pa\be}P_H^{(\al,\be)}$, $\frac\pa{\pa\al}Q_H^{(\al,\be)}$,
$\frac\pa{\pa\be}Q_H^{(\al,\be)}$, respectively, at $(\al,\be)=(0,0)$.
\begin{proposition}
For $m>n$ the matrix entries of $A_Q,B_Q,A_P,B_P$ as occurring in \eqref{55}
are explicitly given by
\begin{align}
(A_Q)_{m,n}&=-\,\frac1{m-n}\,\frac{(-1-x)^{m-n}}{2^{m-n}}\,,\quad
(B_Q)_{m,n}=-\,\frac1{m-n}\,\frac{(1-x)^{m-n}}{2^{m-n}},
\label{54}\\
(A_P)_{m,n}&=\frac1{m-n}\,P_{m-n}^{(0,-1)}(x),\qquad\quad\;\;
(B_P)_{m,n}=\frac1{m-n}\,P_{m-n}^{(-1,0)}(x).
\label{56}
\end{align}
\end{proposition}
\Proof
First note that by \eqref{53}, \eqref{52}, \eqref{57} and \eqref{49} we have
\[
(A_Q)_{m,n}(x)=(-1)^{m-n}(B_Q)_{m,n}(-x),\qquad
(A_P)_{m,n}(x)=(-1)^{m-n}(B_P)_{m,n}(-x).
\]
Thus for \eqref{54} we only have to compute $B_Q$.
We get from \eqref{10} that, for $m>n$,
\[
P_{m-n}^{(n-m,\be+n-m)}(x)=
\frac{(\be+n-m+1)_{m-n}}{(m-n)!}\,\frac{(x-1)^{m-n}}{2^{m-n}}\,.
\]
Differentiation with respect to $\be$ and putting $\be=0$ yields
$(B_Q)_{m,n}$ and, by \eqref{49} also $(A_Q)_{m,n}$, as given in \eqref{54}.

For \eqref{56} we only have to compute $A_P$. Denote the two equal sides
of the generating function \eqref{50} by $f^{(\al,\be)}(w)$. Then
\[
\log(f^{(\al,0)}(w))=\al \log\left(\frac2{1-w+R}\right)-\log R,\qquad
\frac\pa{\pa\al}\,\log(f^{(\al,0)}(w))=
\log\left(\frac2{1-w+R}\right),
\]
\begin{align*}
\frac{\pa^2}{\pa w\,\pa\al}\,\log(f^{(\al,0)}(w))&=\frac{R-w+x}{R(R-w+1)}=
\frac{R+w+1}{2wR}-\frac1w\\
&=\frac1w\,(f^{(0,-1)}(w)-1)
=\sum_{n=1}^\iy P_n^{(0,-1)}(x)\,w^{n-1}.
\end{align*}
Since
$\frac\pa{\pa\al}\,\log(f^{(\al,0)}(0))=0$, we conclude that
\begin{align*}
\frac\pa{\pa\al}\,\log(f^{(\al,0)}(w))&=
\sum_{n=1}^\iy n^{-1} P_n^{(0,-1)}(x)\,w^n,\\
\frac\pa{\pa\al}f^{(\al,0)}(w)\Big|_{\al=0}&=
\sum_{n=0}^\iy\left(\sum_{k=1}^n k^{-1}
P_k^{(0,-1)}(x)\,P_{n-k}(x)\right)w^n,\\
\frac\pa{\pa\al}P_n^{(\al,0)}(x)\Big|_{\al=0}&=
\sum_{k=1}^n k^{-1} P_k^{(0,-1)}(x)\,P_{n-k}(x).
\end{align*}
Since
\[
P_{m-n}^{(\al,0)}(x)=P_{m,n}^{(\al,0)}=\sum_{k=n}^m(P_H^{(\al,0)})_{m,k}\,
P_{k-n}(x),
\]
we have
\[
\frac\pa{\pa\al}P_{m-n}^{(\al,0)}(x)\Big|_{\al=0}=
\sum_{k=n}^{m-1}(A_P)_{m,k}\,P_{k-n}(x)
\]
We conclude that $(A_P)_{m,n}$ is as given by \eqref{56}.\qed
\bPP
Compare the definition \eqref{20} of $L_{m,n}^{(\al,\be)}$ with the definitions
\eqref{52} of $P^{(\al,\be)}_{m,n}$ and $Q^{(\al,\be)}_{m,n}$. It follows that
\begin{equation}
Q_{m,n}^{(\al+m,\be+m)}=L_{m,n}^{(\al,\be)}=P_{m,n}^{(\al+n,\be+n)}.
\label{62}
\end{equation}
\begin{proposition}
We have
\begin{equation}
P_H^{(\al,\be)} L^{(0,0)}=L^{(\al,\be)}=L^{(0,0)} Q_H^{(\al,\be)}.
\label{58}
\end{equation}
\end{proposition}
\Proof
The second equality follows from
\[
L_{m,n}^{(\al,\be)}=(Q_H^{(\al+m,\be+m)})_{m,n}=
\sum_{k=n}^m (Q_H^{(m,m)})_{m,k} (Q_H^{(\al,\be)})_{k,n}=
\sum_{k=n}^m L_{m,k}^{(0,0)} (Q_H^{(\al,\be)})_{k,n}.
\]
The first equality follows from
\begin{align*}
(P^{(0,0)}P_H^{(\al,\be)} L^{(0,0)})_{m,n}&=(P^{(\al,\be)}L^{(0,0)})_{m,n}
=\sum_{k=n}^m P_{m,k}^{(\al,\be)} L_{k,n}^{(0,0)}
=\sum_{k=n}^m P_{m,k}^{(\al,\be)} P_{k,n}^{(n,n)}\\
&=\sum_{j=n}^m P_{m,j}^{(0,0)} P_{j,n}^{(\al+n,\be+n)}=
\sum_{j=n}^m P_{m,j}^{(0,0)} L_{j,n}^{(\al,\be)}
=(P^{(0,0)} L^{(\al,\be)})_{m,n}.
\end{align*}
Here we used \eqref{53}, \eqref{46} and \eqref{62}.\qed
\begin{remark}
It follows from \eqref{58} that
\begin{equation*}
P_H^{(\al,\be)}=L^{(0,0)} Q_H^{(\al,\be)} (L^{(0,0)})^{-1}.
\end{equation*}
So the two-parameter groups $(\al,\be)\mapsto  P_H^{(\al,\be)}$ and
$(\al,\be)\mapsto  Q_H^{(\al,\be)}$ are conjugate by $L^{(0,0)}$ in the group
$\FSL_\iy$.
Note that the inverse of $L^{(0,0)}$ is given by \eqref{20} as
\[
M_{m,n}^{(0,0)}=\frac mn\,P_{m-n}^{(-m,-m)}(x).
\]
\end{remark}
\section{Biorthogonal systems with respect to bilinear forms}
\label{81}
In this section we build on the results of Theorem \ref{18},
Remark \ref{78} and Section \ref{60} in order to obtain systems of
functions on $\ZZ$, involving the functions \eqref{75}, which are biorthogonal
with respect to some explicit bilinear form on $\ZZ$.

If $\al_1+\al_2=-n=\be_1+\be_2$ then the \LHS\ of \eqref{45}
can be evaluated by an elementary expression, where we will use
Proposition \ref{65} and formula \eqref{68}.
Indeed, if $n=2m>0$ then
\begin{multline}
\sum_{k=0}^{2m} P_{2m-k}^{(\al_1,\be_1)}(x)\,P_k^{(\al_2,\be_2)}(x)=
\sum_{k=0}^{2m} P_{2m-k}^{(-m,-m)}(x)\,P_k^{(-m,-m)}(x)\\
=2P_{2m}^{(-m,-m)}(x)=2\left(\frac{x^2-1}4\right)^m\qquad
(\al_1+\al_2=-2m=\be_1+\be_2),
\label{72}
\end{multline}
and if $n=2m-1$ then
\begin{multline}
\sum_{k=0}^{2m-1} P_{2m-1-k}^{(\al_1,\be_1)}(x)\,P_k^{(\al_2,\be_2)}(x)
=\sum_{k=0}^{2m-1} P_{2m-1-k}^{(-m+1,-m+1)}(x)\,P_k^{(-m,-m)}(x)
=P_{2m}^{(-m+1,-m+1)}(x)
\\+P_{2m-2}^{(-m+1,-m+1)}P_1^{(-m,-m)}(x)
=x\left(\frac{x^2-1}4\right)^{m-1}\quad(\al_1+\al_2=-2m+1=\be_1+\be_2).
\label{73}
\end{multline}
These results can be rephrased as identities in $\FSL_\iy$. Define
$R\in\FSL_\iy$ by
\begin{equation*}
R_{m,n}:=
\begin{cases}
1&\mbox{if $m=n$},\\
2\left(\frac{x^2-1}4\right)^{\half(m-n)}&\mbox{if $m>n$ and $m-n$ even},\\
x\left(\frac{x^2-1}4\right)^{\half(m-n-1)}&\mbox{if $m>n$ and $m-n$ odd}.
\end{cases}
\end{equation*}
Also define $\FSJ\colon\FSL_\iy\to\FSL_\iy$ by
$(\FSJ A)_{m,n}:=A_{-n,-m}$.
Then $\FSJ(AB)=(\FSJ B)(\FSJ A)$.
By \eqref{20} we have for $m>n$ that
\[
\big((\FSJ L^{(-\al,-\be)})L^{(\al,\be)}\big)_{m,n}
=\sum_{k=0}^{m-n} P_{m-n-k}^{(-\al-m,-\be-m)}(x)\,P_k^{(\al+n,\be+n)}(x),
\]
which can be evaluated by \eqref{72} and \eqref{73}.
Thus we have obtained that
\begin{equation}
(\FSJ L^{(-\al,-\be)})\,L^{(\al,\be)}=R.
\label{74}
\end{equation}

\begin{proposition}
The inverse $S$ of $R$ in $\FSL_\iy$
(for which we will also use a notation $\mu_x$) is explicitly given by
\begin{equation}
S_{m,n}=\mu_x(m,n)=\begin{cases}
1&\mbox{if $m=n$},\\
\left(\frac{-1-x}2\right)^{m-n}+\left(\frac{1-x}2\right)^{m-n}&
\mbox{if $m>n$}.
\end{cases}
\label{85}
\end{equation}
\end{proposition}
\Proof
It is sufficient to show that
\[
\sum_{k=0}^nR_{n-k,0}\,S_{k,0}=\de_{n,0}\,.
\]
This follows because the generating functions
\begin{align*}
\sum_{n=0}^\iy R_{n,0}\,w^n&=
\frac{\big(1+\thalf(x+1)w\big)\big(1+\thalf(x-1)w\big)}
{1-\tfrac14(x^2-1)w^2}\,,\\
\sum_{n=0}^\iy S_{n,0}\,w^n&=
\frac{1-\tfrac14(x^2-1)w^2}{\big(1+\thalf(x+1)w\big)\big(1+\thalf(x-1)w\big)}
\end{align*}
are inverse to each other. These generating functions, convergent for
$x\in[-1,1]$, $|w|<1$, are immediately computed by geometric series.\qed
\bPP
From \eqref{74} and $S=R^{-1}$ we obtain that
\begin{equation}
L^{(\al,\be)}S(\FSJ L^{(-\al,-\be)})=I.
\label{77}
\end{equation}
Here we used that in $\FSL_\iy$ the implication
$AB=I\;\Rightarrow\;BA=I$ holds.
Formula \eqref{77}  can be rewritten as
\begin{equation*}
\sum_{k,\ell=-\iy}^\iy
L_{m,k}^{(\al,\be)}\,S_{k,\ell}\,L_{-n,-\ell}^{(-\al,-\be)}=\de_{m,n},
\end{equation*}
where the sum only runs over $k,\ell$ such that $n\le\ell\le k\le m$.
With the notation \eqref{75} and with $\mu_x$ given by \eqref{85}
we have obtained:
\begin{proposition}
\begin{equation}
\sum_{k,\ell=-\iy}^\iy
\phi_m^{(\al,\be,x)}(k)\,\phi_{-n}^{(-\al,-\be,x)}(-\ell)\,
\mu_x(k,\ell)=\de_{m,n},
\label{76}
\end{equation}
where the sum only runs over $k,\ell$ such that $n\le\ell\le k\le m$.
\end{proposition}

It is of interest to compare \eqref{76} with the biorthogonality
relation \eqref{69}.
Formula \eqref{76} can also be considered as a biorthogonality
relation, but this
time with respect to the bilinear form $\mu_x$ on $\ZZ$.

\begin{remark}
From Theorem \ref{18} and formula \eqref{77} we obtain
\begin{equation*}
M_{m,n}^{(\al,\be)}=S(\FSJ L^{(-\al,-\be)}).
\end{equation*}
Equivalently, we obtain from  \eqref{69} and
\eqref{76} that
\begin{equation*}
\psi_{-n}^{(\al,\be,x)}(-k)
=\sum_{\ell\in\ZZ}\mu_x(k,\ell)\,\phi_{-n}^{(-\al,-\be,x)}(-\ell)
\end{equation*}
with sum running over $n\le \ell\le k$.
\end{remark}

If we consider the \LHS\ of \eqref{74} with the two factors interchanged
then we can evaluate it by an earlier result. Indeed,
\[
\big(L^{(\al,\be)}\,(\FSJ L^{(-\al,-\be)})\big)_{m,n}
=\sum_{k=n}^m P_{m-k}^{(\al+k,\be+k)}(x)\,
P_{k-n}^{(-\al-k,-\be-k)}(x)=P_{m-n}(x)
\]
by \eqref{12}. Hence
\begin{equation}
L^{(\al,\be)}\,(\FSJ L^{(-\al,-\be)})=P^{(0,0)}.
\label{79}
\end{equation}
$P^{(0,0)}$ has an inverse $T=(P^{(0,0)})^{-1}$ in $\FSL_\iy$, which was already
computed after \eqref{57} and which we also write as $\nu_x$\,:
\begin{equation}
T_{m,n}=\nu_x(m,n)=\begin{cases}
1&\mbox{if $m=n$},\\
-x&\mbox{if $m=n+1$},\\
\frac{1-x^2}{2(n-1)}\,P_{n-2}^{(1,1)}(x)&\mbox{if $m\ge n+2$}.
\end{cases}
\label{84}
\end{equation}
From \eqref{79} we obtain
\begin{equation}
(\FSJ L^{(-\al,-\be)})\,T\,L^{(\al,\be)}=I.
\label{82}
\end{equation}
With the notation \eqref{75} and with $\nu_x$ given by \eqref{84},
the identity \eqref{82} takes the form
\begin{equation}
\sum_{k,\ell=-\iy}^\iy \phi_{-k}^{(-\al,-\be,x)}(-m)\,
\phi_\ell^{(\al,\be,x)}(n)\,
\nu_x(k,\ell)=\de_{m,n}.
\label{83}
\end{equation}
Just as \eqref{76}, we can consider \eqref{83}
as a biorthogonality relation
for two systems of functions on~$\ZZ$ (the duals of the ones in \eqref{73})
with respect to a bilinear form on $\ZZ$, here $\nu_x$.
\section{Limits of a connection formula for Askey-Wilson polynomials}
\label{35}
Askey-Wilson polynomials \cite{5} are defined by
\begin{equation}
p_n(\cos\tha;a_1,a_2,a_3,a_4\mid q):=\frac{(a_1a_2,a_1a_3,a_1a_4;q)_n}{a_1^n}\,
\qhyp43{q^{-n},a_1a_2a_3a_4q^{n-1},a_1e^{i\tha},a_1e^{-i\tha}}
{a_1a_2,a_1a_3,a_1a_4}{q,q}.
\label{23}
\end{equation}
They are symmetric in $a_1,a_2,a_3,a_4$.
The connection coefficients $c_{n,k}$ in
\begin{equation}
p_n(\cos\tha;b_1,b_2,b_3,a_4\mid q)=\sum_{k=0}^n
c_{n,k}(b_1,b_2,b_3,a_4;a_1,a_2,a_3,a_4\mid q)\,
p_k(\cos\tha;a_1,a_2,a_3,a_4\mid q)
\label{21}
\end{equation}
are explicitly given in Askey \& Wilson \cite[(6.5)]{5}:
\begin{multline}
c_{n,k}(b_1,b_2,b_3,a_4;a_1,a_2,a_3,a_4\mid q)
=\frac{q^{k(k-n)}(q;q)_n}{a_4^{n-k}(q;q)_{n-k}(q;q)_k}\,
\frac{(b_1b_2b_3a_4q^{n-1};q)_k}{(a_1a_2a_3a_4q^{k-1};q)_k}\sLP
\times(b_1a_4q^k,b_2a_4q^k,b_3a_4q^k;q)_{n-k}\,
\qhyp54{q^{k-n},b_1b_2b_3a_4q^{n+k-1},a_1a_4q^k,a_2a_4q^k,a_3a_4q^k}
{b_1a_4q^k,b_2a_4q^k,b_3a_4q^k,a_1a_2a_3a_4q^{2k}}{q,q}.
\label{22}
\end{multline}
See also Ismail \& Zhang \cite[Section 3]{6}
and Ismail \cite[\S16.4]{7}, where the connection coefficients are given
more generally for $a_4\ne b_4$. However, note that in \cite[(3.13)]{6}
and \cite[(16.4.3)]{7}
one should read $c_{n,k}({\bf b,a})$ instead of $c_{n,k}({\bf a,b})$.

Now put
\begin{equation}
a_4:=q^{\al+1}/b_1,\quad
b_3:=q^{\be+1}/b_2
\label{24}
\end{equation}
in \eqref{21} and \eqref{22}, and multiply both sides of \eqref{21}
by $1/(q;q)_n$.
By \eqref{23} and \eqref{10} we see that
\begin{align*}
&\lim_{q\to1}\frac1{(q;q)_n}\,
p_n(\cos\tha;b_1,b_2,q^{\be+1}/b_2,q^{\al+1}/b_1\mid q)\sLP
&\quad=\left(\frac{(1-b_1b_2)(b_2-b_1)}{b_1b_2}\right)^n\frac{(\al+1)_n}{n!}\,
\hyp21{-n,n+\al+\be+1}{\al+1}
{\frac{b_2(1-2b_1\cos\tha+b_1^2)}{(1-b_1b_2)(b_2-b_1)}}\sLP
&\quad=\left(\frac{(1-b_1b_2)(b_2-b_1)}{b_1b_2}\right)^n
P_n^{(\al,\be)}
\left(1-2\,\frac{b_2(1-2b_1\cos\tha+b_1^2)}{(1-b_1b_2)(b_2-b_1)}\right).
\end{align*}
By \eqref{21} we also see that
\begin{align*}
&\lim_{q\to1}p_k(\cos\tha;a_1,a_2,a_3,q^{\al+1}/b_1\mid q)\\
&\quad=\left(\frac{(1-a_1a_2)(1-a_1a_3)(b_1-a_1)}{a_1b_1}\right)^k
\hyp10{-k}{-}{\frac{(b_1-a_1a_2a_3)(1-2a_1\cos\tha+a_1^2)}
{(1-a_1a_2)(1-a_1a_3)(b_1-a_1)}}\sLP
&\quad=\left(\frac{(1-a_1a_2)(1-a_1a_3)(b_1-a_1)
-(b_1-a_1a_2a_3)(1-2a_1\cos\tha+a_1^2)}{a_1b_1}\right)^k.
\end{align*}
For the ${}_5\phi_4$ in \eqref{22} we get
\begin{align*}
&\lim_{q\to1}
\frac{(q^{\al+k+1};q)_{n-k}}{(q;q)_{n-k}}\,
\qhyp54{q^{k-n},q^{n+k+\al+\be+1},q^{\al+k+1}a_1/b_1,
q^{\al+k+1}a_2/b_1,q^{\al+k+1}a_3/b_1}
{q^{\al+k+1},q^{\al+k+1}b_2/b_1,q^{\al+\be+k+2}/(b_1b_2),
q^{\al+2k+1}a_1a_2a_3/b_1}{q,q}\sLP
&\quad=\frac{(\al+k+1)_{n-k}}{(n-k)!}\,
\hyp21{-n+k,n+k+\al+\be+1}{\al+k+1}
{\frac{b_2(b_1-a_1)(b_1-a_2)(b_1-a_3)}{(b_1-b_2)(b_1b_2-1)(b_1-a_1a_2a_3)}}\sLP
&\quad=
P_{n-k}^{(\al+k,\be+k)}\left(1-2\,
\frac{b_2(b_1-a_1)(b_1-a_2)(b_1-a_3)}
{(b_1-b_2)(b_1b_2-1)(b_1-a_1a_2a_3)}\right).
\end{align*}
For the other factors in \eqref{22} we get
\begin{align*}
&\lim_{q\to1}
\frac{q^{k(k-n)}}{a_4^{n-k}(q;q)_k}\,
\frac{(q^{n+\al+\be+1};q)_k}{(q^{\al+k}a_1a_2a_3/b_1;q)_k}\,
(q^{\al+k+1}b_2/b_1,q^{\al+\be+k+2}/(b_1b_2);q)_{n-k}\\
&\quad=
\frac{(n+\al+\be+1)_k}{k!}\,
\left(\frac{b_1}{b_1-a_1a_2a_3}\right)^k
\left(\frac{(b_1-b_2)(b_1b_2-1)}{b_1b_2}\right)^{n-k}.
\end{align*}
Also put
\begin{equation}
x=1-2\,\frac{b_2(1-2b_1\cos\tha+b_1^2)}{(1-b_1b_2)(b_2-b_1)},\quad
y=1-2\,
\frac{b_2(b_1-a_1)(b_1-a_2)(b_1-a_3)}
{(b_1-b_2)(b_1b_2-1)(b_1-a_1a_2a_3)}\,.
\label{28}
\end{equation}
Then we obtain the following limit case of \eqref{21} as $q\to1$:
\begin{equation}
P_n^{(\al,\be)}(x)=\sum_{k=0}^n \frac{(n+\al+\be+1)_k}{k!}\,
P_{n-k}^{(\al+k,\be+k)}(y)\,\left(\frac{x-y}2\right)^k.
\label{27}
\end{equation}

Now interchange  the $a$ and $b$ parameters in \eqref{21}:
\begin{equation}
p_n(\cos\tha;a_1,a_2,a_3,a_4\mid q)=\sum_{k=0}^n
c_{n,k}(a_1,a_2,a_3,a_4;b_1,b_2,b_3,a_4\mid q)\,
p_k(\cos\tha;b_1,b_2,b_3,a_4\mid q),
\label{25}
\end{equation}
and use \eqref{22} with the $a$ and $b$ parameters interchanged and with
the order of summation reversion formula \cite[Exercise 1.4(ii)]{8}
applied to the ${}_5\phi_4$:
\begin{align}
&c_{n,k}(a_1,a_2,a_3,a_4;b_1,b_2,b_3,a_4\mid q)
=\frac{(-1)^{n-k}q^{-\half(n-k)(n+k-1)}(q;q)_n}{a_4^{n-k}(q;q)_{n-k}(q;q)_k}
\nonumber\sLP
&\quad\times
\frac{(a_1a_2a_3a_4q^{n-1};q)_n}{(b_1b_2b_3a_4q^{k-1};q)_k
(b_1b_2b_3a_4q^{2k};q)_{n-k}}\,
(b_1a_4q^k,b_2a_4q^k,b_3a_4q^k;q)_{n-k}\nonumber\sLP
&\quad\times
\qhyp54{q^{k-n},q^{1-k-n}/(b_1b_2b_3a_4),q^{1-n}/(a_1a_4),
q^{1-n}/(a_2a_4),q^{1-n}/(a_3a_4)}
{q^{1-n}/(b_1a_4),q^{1-n}/(b_2a_4),q^{1-n}/(b_3a_4),q^{1-2n}/
(a_1a_2a_3a_4)}{q,q}.
\label{26}
\end{align}
Now substitute \eqref{24} in \eqref{25} and \eqref{26} and let $x$ and $y$
be given by
\eqref{28}.
By similar computations as for obtaining \eqref{27} we get as the limit of
\eqref{25} for $q\to1$ the following identity:
\begin{equation}
\left(\frac{x-y}2\right)^n
=\sum_{k=0}^n\frac{\al+\be+2k+1}{\al+\be+k+1}
\frac{n!}{(\al+\be+k+2)_n}\,
P_{n-k}^{(-\al-n-1,-\be-n-1)}(y)\,P_k^{(\al,\be)}(x).
\label{29}
\end{equation}
Formula \eqref{29} was earlier given by
J. Koekoek \& R. Koekoek \cite[(21)]{9}.
As an alternative to their direct derivation (independently of the
Askey-Wilson connection coefficients) one can compute that
\begin{multline*}
\frac{\int_{-1}^1(x-y)^n\,P_k^{(\al,\be)}(x)\,(1-x)^\al(1+x)^\be\,dx}
{\int_{-1}^1(P_k^{(\al,\be)}(x))^{2}\,(1-x)^\al(1+x)^\be\,dx}\\
=\frac{\al+\be+2k+1}{\al+\be+k+1}
\frac{2^n\,n!}{(\al+\be+k+2)_n}\,
P_{n-k}^{(-\al-n-1,-\be-n-1)}(y)
\end{multline*}
by substituting the Rodrigues formula for Jacobi polynomials in the numerator
on the \LHS, then performing repeated integration by parts, then using Euler's
integral representation for hypergeometric functions and finally reversing
the order of summation in the resulting terminating hypergeometric series.
\begin{remark}
Formula \eqref{29} can also be obtained as the special case $\nu=-n$ of
formula (3.1) in Cohl \cite{11}. In that formula he gives an explicit expansion
of $(z-x)^{-\nu}$ in terms of Jacobi polynomials $P_n^{(\al,\be)}(x)$,
where the
expansion coefficients turn out to be constant multiples of the expressions
$(z-1)^{\al+1-\nu}(z+1)^{\be+1-\nu} Q_{n+\nu-1}^{(\al+1-\nu,\be+1-\nu)}(z)$
(the $Q$-functions being Jacobi functions of the second kind).
His case $\nu=1$ occurs in Szeg\H{o} \cite[(9.2.1)]{12}.
Cohl's formula could also have been proved
in the way just sketched for \eqref{29}.
\end{remark}

From \eqref{27} and \eqref{29} we see (as also observed in \cite{9})
that $AB=I=BA$, where $A$ and $B$ are
the lower triangular matrices given for $m\ge n\ge0$ by
\begin{equation}
\label{36}
\begin{split}
A_{m,n}&=\frac{(\al+\be+m+1)_n}{n!}\,P_{m-n}^{(\al+n,\be+n)}(y),\\
B_{m,n}&=\frac{\al+\be+2n+1}{\al+\be+n+1}\,
\frac{m!}{(\al+\be+n+2)_m}\,P_{m-n}^{(-\al-m-1,-\be-m-1)}(y).
\end{split}
\end{equation}
In particular, we obtain from $AB=I$ the identities \eqref{31}
and \eqref{30},
while conversely from \eqref{31} with $(\al,\be)$ running through all
$(\al+j,\be+j)$ ($j\in\Znonneg$)
the full set of scalar identities in $AB=I$ for $(\al,\be)$ can be derived.
Similarly we obtain from $BA=I$ that
\begin{equation}
\sum_{k=0}^n \frac{\al+\be+2k+1}{\al+\be+1}
\frac{(\al+\be+1)_k}{(\al+\be+n+2)_k}\,P_k^{(\al,\be)}(y)\,
P_{n-k}^{(-\al-n-1,-\be-n-1)}(y)=\de_{n,0}.
\label{32}
\end{equation}
Formula \eqref{32} also follows from \eqref{29} by putting $x=y$, as already
observed in \cite[(22)]{9}. Conversely (see
\cite[p.13]{9}),
from \eqref{32} with $(\al,\be)$ running through all
$(\al+j,\be+j)$ ($j\in\Znonneg$)
the full set of scalar identities in $BA=I$ for $(\al,\be)$ can be derived.

\subsection*{Note added in proof}
The application in \cite{2} of the case $\alpha=\beta=\tfrac{1}{2}$ of Theorem \ref{18} has now been generalized to general parameter values $\alpha=\beta>-\tfrac{1}{2}$ in the preprint E. Koelink, A.M. de los Rios, P. Rom\'an, {\em Matrix-valued Gegenbauer polynomials},
{\tt arXiv:1403.2938}.

\quad
\begin{footnotesize}
\begin{quote}
L. Cagliero,
CIEM-CONICET, FAMAF-Universidad Nacional de C\'ordoba, C\'ordoba, Argentina;\\
email: {\tt cagliero@famaf.unc.edu.ar}
\bLP
{ T. H. Koornwinder, Korteweg-de Vries Institute, University of
 Amsterdam,\\
 P.O.\ Box 94248, 1090 GE Amsterdam, The Netherlands;\\
email: }{\tt T.H.Koornwinder@uva.nl}
\end{quote}
\end{footnotesize}

\end{document}